\documentclass[12pt]{article}%
\usepackage[mathlines]{lineno}
\usepackage{amsmath}
\usepackage{amsmath}
\usepackage{amssymb}
\usepackage{dsfont}
\usepackage{cite}
\newsymbol \blackbox 1004
\newcommand{\eh}{\hfill}\newlength{\sperr}

\newenvironment{proof}{{\settowidth{\sperr}{\bf\rm
Proof}%
\par\addvspace{0.3cm}\noindent\parbox[t]{1.3\sperr}
{\bf\rm P\eh r\eh o\eh o\eh f\eh }%
}}{\nopagebreak\mbox{}
$\blackbox$\par\addvspace{0.3cm}}

\def\nn{\nonumber}

\def\bt{\beta}

\def\Lam{\Lambda}
\def\s{\sigma}

\def\ze{\zeta}

\def\vt{\vartheta}

\def\wt{\widetilde}

\def\p{\partial}

\def\BC{{\mathbb C}}

\def\BR{{\mathbb R}}
\def\BN{{\mathbb N}}

\def\cla{{\mathcal A}}

\def\clj{{\mathcal J}}

\def\cli{{\mathcal I}}

\def\diag{\mathrm{diag}}

\newcommand{\I}{\mathrm{i}}

\newtheorem{Pa}{Paper}[section]
\newtheorem{Tm}[Pa]{{\bf Theorem}}

\newtheorem{Rk}[Pa]{{\bf Remark}}

\newtheorem{Ee}[Pa]{{\bf Example}}

\newtheorem{Pn}[Pa]{{\bf Proposition}}

\usepackage{hyperref}

\title{Hamiltonian systems  with several space variables: \\  dressing, explicit solutions and energy relations}

\author{Alexander Sakhnovich \footnote{Faculty of Mathematics,
University
of
Vienna, 
Oskar-Morgenstern-Platz 1, A-1090, Vienna,
Austria.
E-mail address: oleksandr.sakhnovych@univie.ac.at}}

\date{}
\begin{document}
\maketitle

\begin{abstract}  
We construct so-called Darboux transformations and solutions of the dynamical Hamiltonian systems
with several space variables $\frac{\p  \psi}{\p t}=\sum_{k=1}^r  H_k(t)\frac{\p \psi}{\p \ze_k}\,$  $( H_k(t)= H_k(t)^*)$.
 In particular, such systems are analogs  of the  port-Hamiltonian systems in the
important and insufficiently studied case of  several space variables.
The corresponding energy relations are written down.  The method is illustrated by several examples,
where explicit solutions are given. 
\end{abstract}


\vspace{0.2em}


\section{Introduction} \label{intro}
\setcounter{equation}{0}
{\bf 1.} In this paper, we study dynamical
 systems of the form
\begin{align} &       \label{i1}
\frac{\p \psi }{\p t} =\clj \psi, \quad 
 \clj:=\sum_{k=1}^r H_k(t)\frac{\p}{\p \ze_k} \quad \, \big(H_k(t)=H_k(t)^*\big),
\end{align}
where $\psi=\psi(t,\ze_1,\ldots,\ze_r)$ may be a vector or an $m\times n$ matrix, $H_k(t)$ are $m \times m$ matrix-valued functions (matrix functions) and
$H_k(t)^*$ stands for the conjugate transpose of the matrix $H_k(t)$.  The paper may be considered  as an essential development of   
\cite{MST, ALS-IEOT11, ALSarx}. We note that systems $\frac{\p \psi}{\p t}=\clj \psi$ with symmetric operators $\clj$ appear  in \cite{ALSarx}. Here,
we deal with an equally (or more) important case of skew-symmetric $\clj$. 
In particular, system \eqref{i1} with skew-symmetric $\clj$ is of essential interest in the control theory and in the theory of port-Hamiltonian systems.
For instance, the case of the constant matrices $H_k(t)\equiv {\rm const}$
 was dealt with in well-known thesis \cite[Chapter~8]{Vill}.
So-called port-Hamiltonian systems with one space variable have been actively studied
but much less results are known for the case of several space variables.

System \eqref{i1} is of general interest as well. {\it We call sometimes matrix functions $H_k(t)$ Hamiltonians without requiring} $H_k(t)\geq 0$ $($i.e., without requiring that
the eigenvalues of $H_k$ are nonnegative$)$. This terminology comes from the one dimensional canonical systems,
see the references in \cite{GoKr, ALS17, SaSaR, SaL2}.

{\bf 2.}  We construct so-called Darboux transformations and explicit solutions of  the dynamical systems of the form \eqref{i1}.
Explicit solutions are essential in modelling  and optimization problems. Dressing procedure 
(also called B\"acklund-Darboux transformation or commutation method)  is a fruitful tool in the construction
of the explicit solutions of linear and integrable nonlinear differential equations
(see, e.g.,  \cite{Cies3, CoI, Ge, GeT,  Gu, KaG,  KoSaTe, Mar, MS, Mi, SaSaR, ZM} and numerous references therein).
Furthermore, dressing  may be considered as an important part of the general symmetry theory
\cite{Mont, Olv, RP}.
Here, we use 
our GBDT ({\it generalised B\"acklund-Darboux transformation}) approach, which was initially introduced in
\cite{SaA94}.
For further references on GBDT, see, for instance, \cite{Cies3, FKKS, KaS, KoSaTe, ALS17, ALSarx, SaSaR}.

For systems \eqref{i1}, we assume that $t$ is considered on some interval $\cli$ $(0~\in~\cli)$.
Given an initial system \eqref{i1} (or, equivalently, a set of Hamiltonians $H_k(t)$),  GBDT  is determined by a triple of matrices $\{A,S(0), \Pi(0)\}$ 
satisfying
\begin{align} &       \label{c2}
A, S(0)\in \BC^{n\times n}, \quad \Pi(0)\in \BC^{n\times m}, \quad n,m \in \BN, 
\\ &       \label{c3}
AS(0)-S(0)A^*=\I \Pi(0)\Pi(0)^*, \quad S(0)=S(0)^*,
\end{align}
and by a set of scalars $\{c_k\}$:
\begin{align} &       \label{c3+}
c_k\in \BR, \quad c_k\not\in\sigma(A)\quad  (1\leq k\leq r); \quad c_i\not=c_k \,\, {\mathrm{for}} \,\, i\not=k.
\end{align}
Here, as usual, $\BN$ stands for the set of  positive integer numbers, $\BR$ stands for the real numbers,
$\BC$ stands for the complex plane, and $\BC^{n\times m}$ denotes the set of $n\times m $ matrices with complex-valued entries. 
The notation $\I$ stands for the    imaginary unit ($\I^2=-1$) and $\sigma(A)$ stands for the spectrum of the matrix $A$.
The transformed (GBDT-transformed) system has the same structure as the system \eqref{i1}:
\begin{align} &       \label{i2}
\frac{\p \wt \psi}{\p t}=\sum_{k=1}^r \wt H_k(t)\frac{\p \wt \psi}{\p \ze_k}, \quad  \wt H_k(t)= \wt H_k(t)^*.
\end{align}

The construction of the transformed Hamiltonians $\wt H_k(t)$ and $m\times n$ matrix function $\wt \psi$ satisfying \eqref{i2}
is presented in Section \ref{prel}. The corresponding vector solutions of \eqref{i2} have the form $\wt \psi(t,\ze_1,\ldots,\ze_r)h$
$\big(h\in \BC^{n\times 1}\big)$. 

In the cases of canonical systems and port-Hamiltonian systems with one space variable, the energy is usually
introduced in terms of $\int \psi^*H \psi$. However, in the case of several space variables (and taking into account that $H_k$ depend on time) it seems reasonable
to introduce the energy as the standard norm $\big(\int_Vh^*\psi^*\psi h\big)^{1/2}$, $h\in \BC^{n\times 1}$.  
We note that the expressions $\psi^*H_k \psi$ appear in energy relations \eqref{e1} and \eqref{e2} as well. 

Finally, in Section \ref{Ex} we consider some examples and explicit solutions.
Our conclusions are presented in Section \ref{Con}.

{\it Notations.} Many notations have been explained above.  
The notation $I_m$ denotes the $m\times m$ identity matrix
and the matrix inequality $S>0$ ($S\geq 0$)  means  that the eigenvalues of the matrix $S=S^*$ are positive
(nonnegative). In a similar  way, one interprets  the inequalities $S<0$ and $S\leq 0$.
The notation $\diag\{\cla_1, \cla_2, \ldots, \cla_s\}$ means diagonal matrix with the entries (or blocks)
$\cla_1, \cla_2, \ldots$ on the main diagonal.

\section{Preliminaries: GBDT for  systems \eqref{i1}} \label{prel}
\setcounter{equation}{0}
{\bf 1.} 
As explained in the introduction,   GBDT of the initial system \eqref{i1}  (i.e., GBDT for a set of Hamiltonians $H_k(t)$)
is determined by  a triple of matrices $\{A,S(0), \Pi(0)\}$ and a set of scalars $c_k$ ($1\leq k \leq r$).
The transformed Hamiltonians $\{\wt H_k(t)\}$ and solutions of the transformed system \eqref{i2} are expressed
via the generalised eigenfunctions (matrix functions) $\Pi(t)$  and matrix functions $S(t)$, which take values 
$\Pi(0)$ and $S(0)$, respectively, at $t=0$.
In order to   construct $\Pi(t)$ and $S(t)$ and to apply our theory (see \cite[Appendix A]{ALSarx} and references therein),
we write down an auxiliary system
\begin{align}
&       \label{c1}
\frac{dy}{dt}=G(t,z)y, \quad G(t,z):=\I\sum_{k=1}^r(z-c_k)^{-1}H_k(t).
\end{align}
This system  coincides with the system \cite[(19), (22)]{ALSarx} if we put $j=I_m$ in \cite{ALSarx}.
The strong  restrictions on $H_k$, which appear in  \cite{ALSarx}, are non-essential  for the facts
derived below similar to \cite[Section 3]{ALSarx}. Note that spectral theory of the important particular cases of \eqref{c1}
have been studied in \cite{MST, ALS-IEOT11}.

Similar to the case $j=J=I_m$ from \cite{ALSarx}, the matrix functions $\Pi(t)$ and $S(t)$ are determined
(for the case of the initial system \eqref{c1}) by their
values $\Pi(0)$ and $S(0)$ and by the equations
 \begin{align}   &       \label{c5}
 \Pi^{\prime}(t)=-\I\sum_{k=1}^{r} (A -c_k I_n)^{-1}\Pi(t)H_k(t) \quad \Big( \Pi^{\prime}:=\frac{d}{dt}\Pi\Big),
\\  &       \label{c6} 
S^{\prime}(t)=-\sum_{k=1}^{r} (A -c_k I_n)^{-1} \Pi (t)H_k(t)\Pi(t)^*(A^* -c_k I_n)^{-1},
\end{align}
which follow from \cite[(52) and (54)]{ALSarx}.
These relations jointly with \eqref{c3} yield the matrix identity
\begin{align}
&       \label{c7} 
AS(t)-S(t)A^*=\I \Pi(t)\Pi(t)^*,
\end{align}
which also follows from  \cite[(55)]{ALSarx} and \cite[Corollary A.2]{ALSarx}.
In view of the formula \eqref{c1} and \cite[Corollary A.2]{ALSarx}, important
formula \cite[(7.61)]{SaSaR} takes the form
\begin{align} &       \label{c8}
\Big(\Pi(t)^*S(t)^{-1}\Big)^{\prime}=\I\sum_{k=1}^r\wt H_k(t) \Pi(t)^*S(t)^{-1}(A-c_k I_n)^{-1},
 \end{align} 
 where the transformed Hamiltonians $\wt H_k(t)$ are expressed (at each $t$) in terms of the initial
 Hamiltonians $H_k(t)$ and transfer matrix function (in Lev Sakhnovich sense \cite{SaL1, SaL2}) $w_A(t,\cdot)$:
 \begin{align}
&       \label{c9}
w_A(t,z)=I_m-\I \Pi(t)^*S(t)^{-1}(A-z I_n)^{-1}\Pi(t).
\end{align}
Namely, we have (see \cite[(28) and (29)]{ALSarx})
 \begin{align}
&       \label{c10}
\wt H_k(t)=w_A(t,c_k)H_k(t)w_A(t,c_k)^{-1}.
\end{align}
Moreover, according to \cite[(29)]{ALSarx} the matrix functions $w_A(t,c_k)$ are unitary:
 \begin{align}
&       \label{c11}
w_A(t,c_k)^{-1}=w_A(t,c_k)^*.
\end{align}
Here, $w_A(t,z)$ is a so-called Darboux matrix of the auxiliary system \eqref{c1}.
However,  in our study of the dynamical systems of the \eqref{i1} type,
the basic relation is \eqref{c8} (instead of the differential equation on $w_A$ used in the study of systems \eqref{c1}).
Note also that $\wt H_k(t)$ is unitarily similar to $H_k(t)$.

{\bf 2.} Relations \eqref{c8}, \eqref{c10}, and \eqref{c11} yield our next theorem.
 \begin{Tm} \label{TmSV} Let a set of numbers $c_k$ and  a set of  $m\times m$ matrix functions $H_k(t)=H_k(t)^*$ $(1\leq k \leq r)$ as well as a triple
$\{A, \, S(0),\, \Pi(0)\}$ be given. Assume that relations \eqref{c2}--\eqref{c3+} hold
and that $\Pi(t)$ and $S(t)$ satisfy \eqref{c5} and \eqref{c6}, respectively. 

Then, in the points of invertibility of $S(t)$, the matrix function
 \begin{align} &       \label{c12}
\wt \psi(t,\ze_1,\ldots,\ze_r)=\Pi(t)^*S(t)^{-1}\exp\Big\{\I\sum_{k=1}^r\ze_k(A-c_k I_n)^{-1}\Big\}
 \end{align} 
 satisfies the transformed dynamical system
 \begin{align} &       \label{c13}
\frac{\p \wt \psi}{\p t}=\sum_{k=1}^r \wt H_k(t)\frac{\p \wt \psi}{\p \ze_k}, \quad  \wt H_k(t)= \wt H_k(t)^*=w_A(t,c_k)H_k(t)w_A(t,c_k)^*.
\end{align}
 \end{Tm}
 
 The first equality in \eqref{c13}  means that each column of $\wt \psi$ satisfies the same linear dynamical system above.
\begin{Rk} 
Clearly, \eqref{c13} implies that in the case $H_k(t)\geq 0$ we have $\wt H_k(t) \geq 0$ as well.
Moreover, if 
 \begin{align} &       \label{c14}
H_k(t)\geq 0 \quad (1\leq k\leq r),
\end{align}
then, according to \eqref{c6}, we have $S^{\prime}(t)\leq 0$. Thus, if \eqref{c14} holds, it follows that
$S(t)<0$ in the case 
 $$S(0)<0,  \quad \cli=[0,a] \quad (a>0),$$
 and $S(t)>0$ in the case 
 $$S(0)>0, \quad \cli=[-a,0] \quad (a>0).$$
In both cases, $S(t)$ is invertible. 
\end{Rk}
The next proposition may be considered as some kind of conservation law.
\begin{Pn}\label{PnCL} Let relations \eqref{c2}, \eqref{c3}, \eqref{c5}, and \eqref{c6} hold. Then,
in the points of invertibility of $S(t)$ we have
\begin{align} &       \label{d5}
\big(\Pi(t)^*S(t)^{-1}\Pi(t)\big)^{\prime}=\sum_{k=1}^r\big(\wt H_k(t)-H_k(t)\big),
 \end{align}  
where $\wt H_k$ is given by  \eqref{c10}. 
\end{Pn}
\begin{proof}.
Relations \eqref{c8} and \eqref{c9} yield
\begin{align} &       \label{d6}
\big(\Pi(t)^*S(t)^{-1}\big)^{\prime}\Pi(t)=\sum_{k=1}^r\wt H_k(t)\big(I_m-w_A(t,c_k)\big).
 \end{align}  
From \eqref{c5} and \eqref{d6},  it follows that
\begin{align}        \nn
\big(\Pi(t)^*S(t)^{-1}\Pi(t)\big)^{\prime}=&\sum_{k=1}^r\wt H_k(t)\big(I_m-w_A(t,c_k)\big)
\\ &
\label{d7}
+\sum_{k=1}^r\big(w_A(t,c_k)-I_m\big)H_k(t).
 \end{align}   
 Finally, the equalities \eqref{c10} and \eqref{d7} imply \eqref{d5}.
\end{proof}
\section{Energy relations} \label{energy}
\setcounter{equation}{0}
Since $H_k=H_k^*$, one can easily see that the following
energy relation holds for $\psi$ satisfying \eqref{i1}:
 \begin{align} &       \label{e1}
\frac{\p}{\p t}\displaystyle{\int_V} \psi(t,\ze)^*\psi(t,\ze)d\ze
=
\sum_{k=1}^r \displaystyle{\int_V} \frac{\p}{\p \ze_k}\big(\psi(t,\ze)^*H_k(t)\psi(t,\ze)\big)d\ze,
\end{align}
where $V$ is some bounded (for simplicity) domain and $\ze=(\ze_1,\ldots , \ze_r)$, $d\ze$ stands for $d\ze_1 \ldots d\ze_r$.
If $V$ has the simplest form: $a_k\leq \ze_k\leq b_k$ $(1\leq k\leq r)$ 
formula \eqref{e1} may be rewritten as:
 \begin{align}       \nn
\frac{\p}{\p t}\displaystyle{\int_V} \psi(t,\ze)^*\psi(t,\ze)d\ze
= &
\sum_{k=1}^r \left(\int_{V_k^b}\psi(t,\ze)^*H_k(t)\psi(t,\ze)d\mu_k(\ze) \right.
\\ &       \label{e2} \left.
-\int_{V_k^a}\psi(t,\ze)^*H_k(t)\psi(t,\ze)d\mu_k(\ze)\right);
\end{align}
where $V_k^b$ and $V_k^a$ belong to boundary of $V$, more precisely, 
$$V_k^b=V\cap \{\ze: \,\, \ze_k=b_k\},  \quad V_k^a=V\cap \{\ze: \,\, \ze_k=a_k\},$$ 
and $d\mu_k(\ze)$ is obtained by the removing of $d\ze_k$ from the expression $d\ze_1 \ldots d\ze_r$.
In particular, if the right-hand side of \eqref{e2} tends to zero when  $a_k\to -\infty,$ $b_k\to +\infty$ $(1\leq k\leq r)$
formula \eqref{e2} yields
 \begin{align} &    \label{e3}
\frac{\p}{\p t}\displaystyle{\int_{\BR^r}} \psi(t,\ze_1,\ldots,\ze_r)^*\psi(t,\ze_1,\ldots,\ze_r)d\ze
=0.
\end{align}
\begin{Rk}\label{ERtilde} Clearly, if the conditions of  Theorem \ref{TmSV} are fulfilled, one can substitute $\wt\psi$ and $\wt H_k$ $($instead of $\psi$ and $ H_k$, respectively$)$
into \eqref{e1}--\eqref{e3}.
\end{Rk}
Let us (in addition to Remark \ref{ERtilde}) simplify the expression for $\wt\psi^*\wt H_k\wt \psi$.
\begin{Pn}\label{PnGBDT} Let the conditions of  Theorem \ref{TmSV} be fulfilled. Then, we have
 \begin{align} &    \nn
 \wt \psi(t,\xi_1,\ldots,\xi_r)^*\wt H_k(t)\wt \psi(t,\xi_1,\ldots,\xi_r)
 \\  & \nn
 =\exp\Big\{\I\sum_{k=1}^r\xi_k(A-c_k I_n)^{-1}\Big\}^*(A^*-c_k I_n)S(t)^{-1}(A-c_kI_n)\Pi(t)H_k(t)\Pi(t)^*
\\  &    \label{e4} \quad
\times (A^*-c_k I_n)S(t)^{-1}(A-c_kI_n)\exp\Big\{\I\sum_{k=1}^r\xi_k(A-c_k I_n)^{-1}\Big\},
\end{align}
where $\wt \psi$ is given by \eqref{c12} and $\wt H_k$ is given in \eqref{c13}.
\end{Pn}
\begin{proof}. According to \eqref{c12} and \eqref{c13} we have
\begin{align} &    \nn
 \wt \psi(t,\xi_1,\ldots,\xi_r)^*\wt H_k(t)\wt \psi(t,\xi_1,\ldots,\xi_r)
 \\  & \nn
 =\exp\Big\{\I\sum_{k=1}^r\xi_k(A-c_k I_n)^{-1}\Big\}^*S(t)^{-1}\Pi(t)w_A(t,c_k)H_k(t)w_A(t,c_k)^*\Pi(t)^*S(t)^{-1}
\\  &    \label{e5} \quad
\times \exp\Big\{\I\sum_{k=1}^r\xi_k(A-c_k I_n)^{-1}\Big\}.
\end{align}
It follows from \eqref{c9} that
 \begin{align}     \label{e6}
S(t)^{-1}\Pi(t)w_A(t,c_k)= & S(t)^{-1}\Pi(t)
\\ & \nn
-\big(\I S(t)^{-1}\Pi(t)\Pi(t)^*S(t)^{-1}\big)(A-c_kI_n)^{-1}\Pi(t).
\end{align}
Moreover, \eqref{c7} yields
 \begin{align} &    \label{e7}
\I S(t)^{-1}\Pi(t)\Pi(t)^*S(t)^{-1}=S(t)^{-1}(A-c_kI_n)-(A^*-c_kI_n)S(t)^{-1}.
 \end{align} 
 Substituting \eqref{e7} into \eqref{e6}, we obtain
\begin{align} &    \label{e8}
S(t)^{-1}\Pi(t)w_A(t,c_k)=(A^*-c_k I_n)S(t)^{-1}(A-c_kI_n)\Pi(t).
\end{align}
Equalities \eqref{e5} and \eqref{e8} imply \eqref{e4}.
\end{proof}
\section{Examples}\label{Ex}
\setcounter{equation}{0}
Relations \eqref{c9}, \eqref{c12} and \eqref{c13} show that the explicit construction of $\Pi(t)$ and $S(t)$
is the main part of the construction of the transformed dynamical systems \eqref{c13} and their
solutions. That is, given $\Pi(t)$ and $S(t)$, the matrix functions $\wt H_k$ and $\wt \psi$  are standardly
calculated explicitly. Therefore, this section is dedicated to the construction of $\Pi(t)$ and $S(t)$. For this
purpose, we modify for the skew-symmetric $\clj$ the case of \cite[Example 3.4]{ALSarx} (see Example  \ref{Ee1})
as well as consider a  more complicated Example \ref{Ee2}.

\begin{Ee}\label{Ee1}  Consider a simple example of the constant initial Hamiltonians:
\begin{align} &       \label{f1} 
H_k\equiv j=\diag\{I_{m_1}, -I_{m_2}\} \,\, (1\leq k\leq r); \,\, m_1,m_2\in \BN, \,\, m_1+m_2=m.
  \end{align} 
 \end{Ee}
 Here,  we partition $\Pi(0)$ into the $n\times m_1$ and $n\times m_2$ blocks $\vt_1$ and $\vt_2$, respectively: 
 $\Pi(0)=\begin{bmatrix} \vt_1 & \vt_2\end{bmatrix}$. Hence, it  follows from \eqref{c5} that 
 \begin{align} &       \label{f2}
\Pi(t)=\begin{bmatrix}\exp\Big\{-\I \, t\sum_{k=1}^{r} (A -c_k I_n)^{-1}\Big\}\vt_1 & \exp\Big\{\I \, t\sum_{k=1}^{r} (A -c_k I_n)^{-1}\Big\}\vt_2\end{bmatrix}.
 \end{align} 
 If $\s(A)\cap\s(A^*)=\emptyset$, there are unique matrices $C_i$ $(i=1,2)$ such that
 \begin{align} &       \label{f3} 
  AC_i-C_iA^*=\I \vt_i\vt_i^*.
  \end{align} 
  Moreover, in view of  \eqref{f2}, \eqref{f3} and equality $\s(A)\cap\s(A^*)=\emptyset$,  there is a unique solution
  $S(t)$ of \eqref{c7} which is  given by
   \begin{align}  \nn     
  S(t)=&\exp\Big\{-\I \, t\sum_{k=1}^{r} (A -c_k I_n)^{-1}\Big\}C_1\exp\Big\{\I \, t\sum_{k=1}^{r} (A^* -c_k I_n)^{-1}\Big\}
 \\ &  \label{f4} 
  +\exp\Big\{\I \, t\sum_{k=1}^{r} (A -c_k I_n)^{-1}\Big\}C_2\exp\Big\{-\I \, t\sum_{k=1}^{r} (A^* -c_k I_n)^{-1}\Big\}.
  \end{align}
  Equality \eqref{c7} for $S(t)$ of the form \eqref{f4} is checked by direct substitution.  In particular,
  we have $S(0)=C_1+C_2$.

\begin{Ee}\label{Ee2}  Next, assume that
 \begin{align} &       \label{f5} 
r=m, \,\,   H_k\equiv \bt_k^*\bt_k,  \,\, \bt_k  \bt_k^*=1 \,\, (1\leq k\leq m), \,\, \bt_i  \bt_k^*=0 \,\, (i\not=k), 
\end{align} 
that is, $\bt_k$ are constant orthonormal $ 1\times m$ row vectors.
\end{Ee}
Equalities \eqref{f5} yield that the matrix $\bt^*:=\begin{bmatrix} \bt_1^* & \bt_2^* & \ldots & \bt_m^*\end{bmatrix}$ is unitary,
and it suffices to recover $\Pi(t) \bt^*$ instead of  $\Pi(t)$. Using \eqref{f5}, we rewrite \eqref{c5} in the form
 \begin{align} &       \nn
 \Pi(t)^{\prime} \bt^*=-\I \begin{bmatrix}(A-c_1 I_n)^{-1}\Pi(t)\bt_1^*&\, \ldots &\, (A-c_m I_n)^{-1}\Pi(t)\bt_m^*\end{bmatrix}.
  \end{align} 
Thus, we obtain 
\begin{align} &       \label{f6} 
 \Pi(t) \bt^*=\begin{bmatrix} \Lam_1(t)  & \ldots & \Lam_m(t)\end{bmatrix}, \quad \Lam_k(t)=\exp\big\{-\I \, t (A-c_k I_n)^{-1}\big\}\Pi(0)\bt_k^*,
  \end{align} 
where $\Lam_k(t)$ is the $k$-th column  of $\Pi(t)\bt^*$.

Let us assume again $\s(A)\cap\s(A^*)=\emptyset$.  Recall that  the matrix identities (equations) $AC-CA^*=Q$  have, in this case, unique solutions $C$.
Therefore, relation    \eqref{c7} and the first equality in \eqref{f6} show that 
 \begin{align} &       \label{f7} 
 S(t)=\sum_{k=1}^m S_k(t), \quad AS_k(t)-S_k(t)A^*=\I \Lam_k(t)\Lam_k(t)^*.
  \end{align} 
Moreover, the second equality in \eqref{f6} implies that  the matrix function
 \begin{align} &       \label{f8} 
 S_k(t)=\exp\big\{-\I \, t (A-c_k I_n)^{-1}\big\}C_k\exp\big\{\I \, t (A^*-c_k I_n)^{-1}\big\}
  \end{align} 
satisfies the second equality in \eqref{f7} if
 \begin{align} &       \label{f9} 
  AC_k-C_kA^*=\I \Pi(0)\bt_k^*\bt_k\Pi(0).
  \end{align}
In other words, $S(t)$ is given by the first equality in \eqref{f7} and by relations \eqref{f8} and \eqref{f9}.

\section{Conclusion}\label{Con}
In this paper, we constructed  GBDT-transformed dynamical systems \eqref{i2} and their solutions
including the case of explicit solutions. We studied energy relations for these solutions as well.
The studies of solutions and energy relations are important in itself and may be further applied
in optimisation and control problems. For instance, in analogy with port-Hamiltonian theory,  the triples $\{A,S(0), \Pi(0)\}$ and
the scalars $\{c_k\}$, which determine GBDT, could be considered as the so-called port-variables. 
In that case, the optimisation and control would consist
in choosing both systems and their solutions and we work further on this approach.
We note that Hamiltonian systems as well as energy functionals and control problems  for dynamical systems 
are actively studied in mathematical and applied literature
(see, e.g., some references in the recent publications \cite{Dos, Las, Mog, ALS17, ALSarx}).

{\bf Funding details.}  This work was supported by the Austrian Science Fund (FWF) 
 grant, DOI: 10.55776/Y963.

\end{document}